\documentclass{amsart}
\usepackage{amsfonts}

\newtheorem{theorem}{Theorem}
\theoremstyle{plain}

\newtheorem{corollary}{Corollary}

\newtheorem{problem}{Problem}

\newtheorem{remark}{Remark}

\numberwithin{equation}{section}
\input{tcilatex}

\begin{document}
\title[Inequalities for the Numerical Radius of Linear Operators]{Reverse
Inequalities for the Numerical Radius of Linear Operators in Hilbert Spaces }
\author{S.S. Dragomir}
\address{School of Computer Science and Mathematics\\
Victoria University of Technology\\
PO Box 14428, Melbourne City\\
Victoria 8001, Australia.}
\email{sever.dragomir@vu.edu.au}
\urladdr{http://rgmia.vu.edu.au/dragomir}
\date{29 August, 2005}

\begin{abstract}
Some elementary inequalities providing upper bounds for the difference of
the norm and the numerical radius of a bounded linear operator on Hilbert
spaces under appropriate conditions are given.
\end{abstract}

\maketitle

\section{Introduction}

Let $\left( H;\left\langle \cdot ,\cdot \right\rangle \right) $ be a complex
Hilbert space. The \textit{numerical range} of an operator $T$ is the subset
of the complex numbers $\mathbb{C}$ given by \cite[p. 1]{GR}:%
\begin{equation*}
W\left( T\right) =\left\{ \left\langle Tx,x\right\rangle ,\ x\in H,\
\left\Vert x\right\Vert =1\right\} .
\end{equation*}%
The following properties of $W\left( T\right) $ are immediate:

\begin{enumerate}
\item[(i)] $W\left( \alpha I+\beta T\right) =\alpha +\beta W\left( T\right) $
for $\alpha ,\beta \in \mathbb{C};$

\item[(ii)] $W\left( T^{\ast }\right) =\left\{ \bar{\lambda},\lambda \in
W\left( T\right) \right\} ,$ where $T^{\ast }$ is the \textit{adjoint
operator }of $T;$

\item[(iii)] $W\left( U^{\ast }TU\right) =W\left( T\right) $ for any \textit{%
unitary} operator $U.$
\end{enumerate}

The following classical fact about the geometry of the numerical range \cite[%
p. 4]{GR} may be stated:

\begin{theorem}[Toeplitz-Hausdorff]
The numerical range of an operator is convex.
\end{theorem}

An important use of $W\left( T\right) $ is to bound the \textit{spectrum }$%
\sigma \left( T\right) $ of the operator $T$ \cite[p. 6]{GR}:

\begin{theorem}[Spectral inclusion]
The spectrum of an operator is contained in the closure of its numerical
range.
\end{theorem}

The self-adjoint operators have their spectra bounded sharply by the
numerical range \cite[p. 7]{GR}:

\begin{theorem}
The following statements hold true:

\begin{enumerate}
\item[(i)] $T$ is self-adjoint iff $W\left( T\right) $ is real;

\item[(ii)] If $T$ is self-adjoint and $W\left( T\right) =\left[ m,M\right] $
(the closed interval of real numbers $m,M$), then $\left\Vert T\right\Vert
=\max \left\{ \left\vert m\right\vert ,\left\vert M\right\vert \right\} .$

\item[(iii)] If $W\left( T\right) =\left[ m,M\right] ,$ then $m,M\in \sigma
\left( T\right) .$
\end{enumerate}
\end{theorem}

The \textit{numerical radius} $w\left( T\right) $ of an operator $T$ on $H$
is given by \cite[p. 8]{GR}:%
\begin{equation}
w\left( T\right) =\sup \left\{ \left\vert \lambda \right\vert ,\lambda \in
W\left( T\right) \right\} =\sup \left\{ \left\vert \left\langle
Tx,x\right\rangle \right\vert ,\left\Vert x\right\Vert =1\right\} .
\label{1.1}
\end{equation}%
Obviously, by (\ref{1.1}), for any $x\in H$ one has%
\begin{equation}
\left\vert \left\langle Tx,x\right\rangle \right\vert \leq w\left( T\right)
\left\Vert x\right\Vert ^{2}.  \label{1.2}
\end{equation}

It is well known that $w\left( \cdot \right) $ is a norm on the Banach
algebra $B\left( H\right) $ of all bounded linear operators $T:H\rightarrow
H,$ i.e.,

\begin{enumerate}
\item[(i)] $w\left( T\right) \geq 0$ for any $T\in B\left( H\right) $ and $%
w\left( T\right) =0$ if and only if $T=0;$

\item[(ii)] $w\left( \lambda T\right) =\left\vert \lambda \right\vert
w\left( T\right) $ for any $\lambda \in \mathbb{C}$ and $T\in B\left(
H\right) ;$

\item[(iii)] $w\left( T+V\right) \leq w\left( T\right) +w\left( V\right) $
for any $T,V\in B\left( H\right) .$
\end{enumerate}

This norm is equivalent with the operator norm. In fact, the following more
precise result holds \cite[p. 9]{GR}:

\begin{theorem}[Equivalent norm]
For any $T\in B\left( H\right) $ one has%
\begin{equation}
w\left( T\right) \leq \left\Vert T\right\Vert \leq 2w\left( T\right) .
\label{1.3}
\end{equation}
\end{theorem}

Let us now look at two extreme cases of the inequality (\ref{1.3}). In the
following $r\left( t\right) :=\sup \left\{ \left\vert \lambda \right\vert
,\lambda \in \sigma \left( T\right) \right\} $ will denote the \textit{%
spectral radius} of $T$ and $\sigma _{p}\left( T\right) =\left\{ \lambda \in
\sigma \left( T\right) ,\ Tf=\lambda f\text{ \ for some \ }f\in H\right\} $
the \textit{point spectrum }of $T.$

The following results hold \cite[p.10]{GR}:

\begin{theorem}
\label{t1}We have

\begin{enumerate}
\item[(i)] If $w\left( T\right) =\left\Vert T\right\Vert ,$ then $r\left(
T\right) =\left\Vert T\right\Vert .$

\item[(ii)] If $\lambda \in W\left( T\right) $ and $\left\vert \lambda
\right\vert =\left\Vert T\right\Vert ,$ then $\lambda \in \sigma _{p}\left(
T\right) .$
\end{enumerate}
\end{theorem}

To address the other extreme case $w\left( T\right) =\frac{1}{2}\left\Vert
T\right\Vert ,$ we can state the following sufficient condition in terms of
(see \cite[p. 11]{GR})%
\begin{equation*}
R\left( T\right) :=\left\{ Tf,\ f\in H\right\} \quad \text{and}\quad R\left(
T^{\ast }\right) :=\left\{ T^{\ast }f,\ f\in H\right\} .
\end{equation*}

\begin{theorem}
\label{t.e}If $R\left( T\right) \perp R\left( T^{\ast }\right) ,$ then $%
w\left( T\right) =\frac{1}{2}\left\Vert T\right\Vert .$
\end{theorem}

It is well-known that the two-dimensional shift%
\begin{equation*}
S_{2}=\left[ 
\begin{array}{ll}
0 & 0 \\ 
1 & 0%
\end{array}%
\right] ,
\end{equation*}%
has the property that $w\left( T\right) =\frac{1}{2}\left\Vert T\right\Vert
. $

The following theorem shows that some operators $T$ with $w\left( T\right) =%
\frac{1}{2}\left\Vert T\right\Vert $ have $S_{2}$ as a component \cite[p. 11]%
{GR}:

\begin{theorem}
If $w\left( T\right) =\frac{1}{2}\left\Vert T\right\Vert $ and $T$ attains
its norm, then $T$ has a two-dimensional reducing subspace on which it is
the shift $S_{2}.$
\end{theorem}

For other results on numerical radius, see \cite{H}, Chapter 11.

The main aim of the present paper is to point out some upper bounds for the
nonnegative difference%
\begin{equation*}
\left\Vert T\right\Vert -w\left( T\right) \qquad \left( \left\Vert
T\right\Vert ^{2}-\left( W\left( T\right) \right) ^{2}\right)
\end{equation*}%
under appropriate assumptions for the bounded linear operator $%
T:H\rightarrow H.$

\section{The Results}

The following results may be stated:

\begin{theorem}
\label{t.1}Let $T:H\rightarrow H$ be a bounded linear operator on the
complex Hilbert space $H.$ If $\lambda \in \mathbb{C}\backslash \left\{
0\right\} $ and $r>0$ are such that%
\begin{equation}
\left\Vert T-\lambda I\right\Vert \leq r,  \label{2.1}
\end{equation}%
where $I:H\rightarrow H$ is the identity operator on $H,$ then%
\begin{equation}
\left( 0\leq \right) \left\Vert T\right\Vert -w\left( T\right) \leq \frac{1}{%
2}\cdot \frac{r^{2}}{\left\vert \lambda \right\vert }.  \label{2.2}
\end{equation}
\end{theorem}

\begin{proof}
For $x\in H$ with $\left\Vert x\right\Vert =1,$ we have from (\ref{2.1}) that%
\begin{equation*}
\left\Vert Tx-\lambda x\right\Vert \leq \left\Vert T-\lambda I\right\Vert
\leq r,
\end{equation*}%
giving%
\begin{equation}
\left\Vert Tx\right\Vert ^{2}+\left\vert \lambda \right\vert ^{2}\leq 2\func{%
Re}\left[ \overline{\lambda }\left\langle Tx,x\right\rangle \right]
+r^{2}\leq 2\left\vert \lambda \right\vert \left\vert \left\langle
Tx,x\right\rangle \right\vert +r^{2}.  \label{2.3}
\end{equation}%
Taking the supremum over $x\in H,$ $\left\Vert x\right\Vert =1$ in (\ref{2.3}%
) we get the following inequality that is of interest in itself:%
\begin{equation}
\left\Vert T\right\Vert ^{2}+\left\vert \lambda \right\vert ^{2}\leq
2w\left( T\right) \left\vert \lambda \right\vert +r^{2}.  \label{2.4}
\end{equation}%
Since, obviously,%
\begin{equation}
\left\Vert T\right\Vert ^{2}+\left\vert \lambda \right\vert ^{2}\geq
2\left\Vert T\right\Vert \left\vert \lambda \right\vert ,  \label{2.5}
\end{equation}%
hence by (\ref{2.4}) and (\ref{2.5}) we deduce the desired inequality (\ref%
{2.2}).
\end{proof}

\begin{remark}
If the operator $T:H\rightarrow H$ is such that $R\left( T\right) \perp
R\left( T^{\ast }\right) ,$ $\left\Vert T\right\Vert =1$ and $\left\Vert
T-I\right\Vert \leq 1$, then the equality case holds in (\ref{2.2}). Indeed,
by Theorem \ref{t.e}, we have in this case $w\left( T\right) =\frac{1}{2}%
\left\Vert T\right\Vert =\frac{1}{2}$ and since we can choose in Theorem \ref%
{t.1}, $\lambda =1,$ $r=1,$ then we get in both sides of (\ref{2.2}) the
same quantity $\frac{1}{2}.$
\end{remark}

\begin{problem}
Find the bounded linear operators $T:H\rightarrow H$ with $\left\Vert
T\right\Vert =1,$ $R\left( T\right) \perp R\left( T^{\ast }\right) $ and $%
\left\Vert T-\lambda I\right\Vert \leq \left\vert \lambda \right\vert ^{%
\frac{1}{2}}.$
\end{problem}

The following corollary may be stated:

\begin{corollary}
\label{c.1}Let $A:H\rightarrow H$ be a bounded linear operator and $\varphi
,\phi \in \mathbb{C}$ with $\phi \neq -\varphi ,\varphi .$ If%
\begin{equation}
\func{Re}\left\langle \phi x-Ax,Ax-\varphi x\right\rangle \geq 0\quad \text{%
for any}\quad x\in H,\ \left\Vert x\right\Vert =1  \label{2.6}
\end{equation}%
then%
\begin{equation}
\left( 0\leq \right) \left\Vert A\right\Vert -w\left( A\right) \leq \frac{1}{%
4}\cdot \frac{\left\vert \phi -\varphi \right\vert ^{2}}{\left\vert \phi
+\varphi \right\vert }.  \label{2.7}
\end{equation}
\end{corollary}

\begin{proof}
Utilising the fact that in any Hilbert space the following two statements
are equivalent:

\begin{enumerate}
\item[(i)] $\func{Re}\left\langle Z-x,x-z\right\rangle \geq 0,$ $x,z,Z\in H;$

\item[(ii)] $\left\Vert x-\frac{z+Z}{2}\right\Vert \leq \frac{1}{2}%
\left\Vert Z-z\right\Vert ,$
\end{enumerate}

we deduce that (\ref{2.6}) is equivalent to%
\begin{equation}
\left\Vert Ax-\frac{\phi +\varphi }{2}\cdot Ix\right\Vert \leq \frac{1}{2}%
\left\vert \phi -\varphi \right\vert  \label{2.8}
\end{equation}%
for any $x\in H,$ $\left\Vert x\right\Vert =1,$ which in its turn is
equivalent with the operator norm inequality:%
\begin{equation}
\left\Vert A-\frac{\phi +\varphi }{2}\cdot I\right\Vert \leq \frac{1}{2}%
\left\vert \phi -\varphi \right\vert .  \label{2.9}
\end{equation}%
Now, applying Theorem \ref{t.1} for $T=A,$ $\lambda =\frac{\varphi +\phi }{2}
$ and $r=\frac{1}{2}\left\vert \Gamma -\gamma \right\vert ,$ we deduce the
desired result (\ref{2.7}).
\end{proof}

\begin{remark}
Following \cite[p. 25]{GR}, we say that an operator $B:H\rightarrow H$ is
accreative, if $\func{Re}\left\langle Bx,x\right\rangle \geq 0$ for any $%
x\in H.$ One may observe that the assumption (\ref{2.6}) above is then
equivalent with the fact that the operator $\left( A^{\ast }-\bar{\varphi}%
I\right) \left( \phi I-A\right) $ is accreative.
\end{remark}

Perhaps a more convenient sufficient condition in terms of positive
operators is the following one:

\begin{corollary}
\label{c.2}Let $\varphi ,\phi \in \mathbb{C}$ with $\phi \neq -\varphi
,\varphi $ and $A:H\rightarrow H$ a bounded linear operator in $H.$ If $%
\left( A^{\ast }-\bar{\varphi}I\right) \left( \phi I-A\right) $ is
self-adjoint and%
\begin{equation}
\left( A^{\ast }-\bar{\varphi}I\right) \left( \phi I-A\right) \geq 0
\label{2.10}
\end{equation}%
in the operator order, then%
\begin{equation}
\left( 0\leq \right) \left\Vert A\right\Vert -w\left( A\right) \leq \frac{1}{%
4}\cdot \frac{\left\vert \phi -\varphi \right\vert ^{2}}{\left\vert \phi
+\varphi \right\vert }.  \label{2.11}
\end{equation}
\end{corollary}

The following result may be stated as well:

\begin{corollary}
\label{c.3}Assume that $T,\lambda ,r$ are as in Theorem \ref{t.1}. If, in
addition, there exists $\rho \geq 0$ such that%
\begin{equation}
\left\vert \left\vert \lambda \right\vert -w\left( T\right) \right\vert \geq
\rho ,  \label{2.12a}
\end{equation}%
then%
\begin{equation}
\left( 0\leq \right) \left\Vert T\right\Vert ^{2}-w^{2}\left( T\right) \leq
r^{2}-\rho ^{2}.  \label{2.13}
\end{equation}
\end{corollary}

\begin{proof}
From (\ref{2.4}) of Theorem \ref{t.1}, we have%
\begin{align}
\left\Vert T\right\Vert ^{2}-w^{2}\left( T\right) & \leq r^{2}-w^{2}\left(
T\right) +2w\left( T\right) \left\vert \lambda \right\vert -\left\vert
\lambda \right\vert ^{2}  \label{2.14} \\
& =r^{2}-\left( \left\vert \lambda \right\vert -w\left( T\right) \right)
^{2}.  \notag
\end{align}%
On utilising (\ref{2.4}) and (\ref{2.12a}) we deduce the desired inequality (%
\ref{2.13}).
\end{proof}

\begin{remark}
In particular, if $\left\Vert T-\lambda I\right\Vert \leq r$ and $\left\vert
\lambda \right\vert =w\left( T\right) ,$ $\lambda \in \mathbb{C}$, then%
\begin{equation}
\left( 0\leq \right) \left\Vert T\right\Vert ^{2}-w^{2}\left( T\right) \leq
r^{2}.  \label{2.15}
\end{equation}
\end{remark}

The following result may be stated as well.

\begin{theorem}
\label{t.2}Let $T:H\rightarrow H$ be a nonzero bounded linear operator on $H$
and $\lambda \in \mathbb{C}\setminus \left\{ 0\right\} ,$ $r>0$ with $%
\left\vert \lambda \right\vert >r.$ If%
\begin{equation}
\left\Vert T-\lambda I\right\Vert \leq r,  \label{3.1}
\end{equation}%
then%
\begin{equation}
\sqrt{1-\frac{r^{2}}{\left\vert \lambda \right\vert ^{2}}}\leq \frac{w\left(
T\right) }{\left\Vert T\right\Vert }\quad \left( \leq 1\right) .  \label{3.2}
\end{equation}
\end{theorem}

\begin{proof}
From (\ref{2.4}) of Theorem \ref{t.1}, we have%
\begin{equation*}
\left\Vert T\right\Vert ^{2}+\left\vert \lambda \right\vert ^{2}-r^{2}\leq
2\left\vert \lambda \right\vert w\left( T\right) ,
\end{equation*}%
which implies, on dividing with $\sqrt{\left\vert \lambda \right\vert
^{2}-r^{2}}>0$ that%
\begin{equation}
\frac{\left\Vert T\right\Vert ^{2}}{\sqrt{\left\vert \lambda \right\vert
^{2}-r^{2}}}+\sqrt{\left\vert \lambda \right\vert ^{2}-r^{2}}\leq \frac{%
2\left\vert \lambda \right\vert w\left( T\right) }{\sqrt{\left\vert \lambda
\right\vert ^{2}-r^{2}}}.  \label{3.3}
\end{equation}%
By the elementary inequality%
\begin{equation}
2\left\Vert T\right\Vert \leq \frac{\left\Vert T\right\Vert ^{2}}{\sqrt{%
\left\vert \lambda \right\vert ^{2}-r^{2}}}+\sqrt{\left\vert \lambda
\right\vert ^{2}-r^{2}}  \label{3.4}
\end{equation}%
and by (\ref{3.3}) we deduce%
\begin{equation*}
\left\Vert T\right\Vert \leq \frac{w\left( T\right) \left\vert \lambda
\right\vert }{\sqrt{\left\vert \lambda \right\vert ^{2}-r^{2}}},
\end{equation*}%
which is equivalent to (\ref{3.2}).
\end{proof}

\begin{remark}
Squaring (\ref{3.2}), we get the inequality%
\begin{equation}
\left( 0\leq \right) \left\Vert T\right\Vert ^{2}-w^{2}\left( T\right) \leq 
\frac{r^{2}}{\left\vert \lambda \right\vert ^{2}}\left\Vert T\right\Vert
^{2}.  \label{3.5}
\end{equation}
\end{remark}

\begin{remark}
Since for any bounded linear operator $T:H\rightarrow H$ we have that $%
w\left( T\right) \geq \frac{1}{2}\left\Vert T\right\Vert ,$ hence (\ref{3.2}%
) would produce a refinement of this classic fact only in the case when%
\begin{equation*}
\frac{1}{2}\leq \left( 1-\frac{r^{2}}{\left\vert \lambda \right\vert ^{2}}%
\right) ^{\frac{1}{2}},
\end{equation*}%
which is equivalent to $r/\left\vert \lambda \right\vert \leq \sqrt{3}/2.$
\end{remark}

The following corollary holds.

\begin{corollary}
\label{c.4}Let $\varphi ,\phi \in \mathbb{C}$ with $\func{Re}\left( \phi 
\bar{\varphi}\right) >0.$ If $T:H\rightarrow H$ is a bounded linear operator
such that either (\ref{2.6}) or (\ref{2.10}) holds true, then:%
\begin{equation}
\frac{2\sqrt{\func{Re}\left( \phi \bar{\varphi}\right) }}{\left\vert \phi
+\varphi \right\vert }\leq \frac{w\left( T\right) }{\left\Vert T\right\Vert }%
\left( \leq 1\right)  \label{3.6}
\end{equation}%
and%
\begin{equation}
\left( 0\leq \right) \left\Vert T\right\Vert ^{2}-w^{2}\left( T\right) \leq
\left\vert \frac{\phi -\varphi }{\phi +\varphi }\right\vert ^{2}\left\Vert
T\right\Vert ^{2}.  \label{3.7}
\end{equation}
\end{corollary}

\begin{proof}
If we consider $\lambda =\frac{\phi +\varphi }{2}$ and $r=\frac{1}{2}%
\left\vert \phi -\varphi \right\vert ,$ then $\left\vert \lambda \right\vert
^{2}-r^{2}=\left\vert \frac{\phi +\varphi }{2}\right\vert ^{2}-\left\vert 
\frac{\phi -\varphi }{2}\right\vert ^{2}=\func{Re}\left( \phi \bar{\varphi}%
\right) >0.$ Now, on applying Theorem \ref{t.2}, we deduce the desired
result.
\end{proof}

\begin{remark}
If $\left\vert \phi -\varphi \right\vert \leq \frac{\sqrt{3}}{2}\left\vert
\phi +\varphi \right\vert ,$ $\func{Re}\left( \phi \bar{\varphi}\right) >0,$
then (\ref{3.6}) is a refinement of the inequality $w\left( T\right) \geq 
\frac{1}{2}\left\Vert T\right\Vert .$
\end{remark}

The following result may be of interest as well.

\begin{theorem}
\label{t.3}Let $T:H\rightarrow H$ be a nonzero bounded linear operator on $H$
and $\lambda \in \mathbb{C}\backslash \left\{ 0\right\} ,$ $r>0$ with $%
\left\vert \lambda \right\vert >r.$ If%
\begin{equation}
\left\Vert T-\lambda I\right\Vert \leq r,  \label{4.1}
\end{equation}%
then%
\begin{equation}
\left( 0\leq \right) \left\Vert T\right\Vert ^{2}-w^{2}\left( T\right) \leq 
\frac{2r^{2}}{\left\vert \lambda \right\vert +\sqrt{\left\vert \lambda
\right\vert ^{2}-r^{2}}}w\left( T\right) .  \label{4.2}
\end{equation}
\end{theorem}

\begin{proof}
From the proof of Theorem \ref{t.1}, we have%
\begin{equation}
\left\Vert Tx\right\Vert ^{2}+\left\vert \lambda \right\vert ^{2}\leq 2\func{%
Re}\left[ \overline{\lambda }\left\langle Tx,x\right\rangle \right] +r^{2}
\label{4.2a}
\end{equation}%
for any $x\in H,$ $\left\Vert x\right\Vert =1.$

If we divide (\ref{4.2a}) by $\left\vert \lambda \right\vert \left\vert
\left\langle Tx,x\right\rangle \right\vert ,$ (which, by (\ref{4.2a}), is
positive) then we obtain%
\begin{equation}
\frac{\left\Vert Tx\right\Vert ^{2}}{\left\vert \lambda \right\vert
\left\vert \left\langle Tx,x\right\rangle \right\vert }\leq \frac{2\func{Re}%
\left[ \overline{\lambda }\left\langle Tx,x\right\rangle \right] }{%
\left\vert \lambda \right\vert \left\vert \left\langle Tx,x\right\rangle
\right\vert }+\frac{r^{2}}{\left\vert \lambda \right\vert \left\vert
\left\langle Tx,x\right\rangle \right\vert }-\frac{\left\vert \lambda
\right\vert }{\left\vert \left\langle Tx,x\right\rangle \right\vert }
\label{4.3}
\end{equation}%
for any $x\in H,$ $\left\Vert x\right\Vert =1.$

If we subtract in (\ref{4.3}) the same quantity $\frac{\left\vert
\left\langle Tx,x\right\rangle \right\vert }{\left\vert \lambda \right\vert }
$ from both sides, then we get%
\begin{align}
& \frac{\left\Vert Tx\right\Vert ^{2}}{\left\vert \lambda \right\vert
\left\vert \left\langle Tx,x\right\rangle \right\vert }-\frac{\left\vert
\left\langle Tx,x\right\rangle \right\vert }{\left\vert \lambda \right\vert }
\label{4.4} \\
& \leq \frac{2\func{Re}\left[ \overline{\lambda }\left\langle
Tx,x\right\rangle \right] }{\left\vert \lambda \right\vert \left\vert
\left\langle Tx,x\right\rangle \right\vert }+\frac{r^{2}}{\left\vert \lambda
\right\vert \left\vert \left\langle Tx,x\right\rangle \right\vert }-\frac{%
\left\vert \left\langle Tx,x\right\rangle \right\vert }{\left\vert \lambda
\right\vert }-\frac{\left\vert \lambda \right\vert }{\left\vert \left\langle
Tx,x\right\rangle \right\vert }  \notag \\
& =\frac{2\func{Re}\left[ \overline{\lambda }\left\langle Tx,x\right\rangle %
\right] }{\left\vert \lambda \right\vert \left\vert \left\langle
Tx,x\right\rangle \right\vert }-\frac{\left\vert \lambda \right\vert
^{2}-r^{2}}{\left\vert \lambda \right\vert \left\vert \left\langle
Tx,x\right\rangle \right\vert }-\frac{\left\vert \left\langle
Tx,x\right\rangle \right\vert }{\left\vert \lambda \right\vert }  \notag \\
& =\frac{2\func{Re}\left[ \overline{\lambda }\left\langle Tx,x\right\rangle %
\right] }{\left\vert \lambda \right\vert \left\vert \left\langle
Tx,x\right\rangle \right\vert }-\left( \frac{\sqrt{\left\vert \lambda
\right\vert ^{2}-r^{2}}}{\sqrt{\left\vert \lambda \right\vert \left\vert
\left\langle Tx,x\right\rangle \right\vert }}-\frac{\sqrt{\left\vert
\left\langle Tx,x\right\rangle \right\vert }}{\sqrt{\left\vert \lambda
\right\vert }}\right) ^{2}-2\frac{\sqrt{\left\vert \lambda \right\vert
^{2}-r^{2}}}{\left\vert \lambda \right\vert }.  \notag
\end{align}%
Since%
\begin{equation*}
\func{Re}\left[ \overline{\lambda }\left\langle Tx,x\right\rangle \right]
\leq \left\vert \lambda \right\vert \left\vert \left\langle
Tx,x\right\rangle \right\vert 
\end{equation*}%
and%
\begin{equation*}
\left( \frac{\sqrt{\left\vert \lambda \right\vert ^{2}-r^{2}}}{\sqrt{%
\left\vert \lambda \right\vert \left\vert \left\langle Tx,x\right\rangle
\right\vert }}-\frac{\sqrt{\left\vert \left\langle Tx,x\right\rangle
\right\vert }}{\sqrt{\left\vert \lambda \right\vert }}\right) ^{2}\geq 0
\end{equation*}%
hence by (\ref{4.4}) we get%
\begin{equation*}
\frac{\left\Vert Tx\right\Vert ^{2}}{\left\vert \lambda \right\vert
\left\vert \left\langle Tx,x\right\rangle \right\vert }-\frac{\left\vert
\left\langle Tx,x\right\rangle \right\vert }{\left\vert \lambda \right\vert }%
\leq \frac{2\left( \left\vert \lambda \right\vert -\sqrt{\left\vert \lambda
\right\vert ^{2}-r^{2}}\right) }{\left\vert \lambda \right\vert }
\end{equation*}%
which gives the inequality%
\begin{equation}
\left\Vert Tx\right\Vert ^{2}\leq \left\vert \left\langle Tx,x\right\rangle
\right\vert ^{2}+2\left\vert \left\langle Tx,x\right\rangle \right\vert
\left( \left\vert \lambda \right\vert -\sqrt{\left\vert \lambda \right\vert
^{2}-r^{2}}\right)   \label{4.5}
\end{equation}%
for any $x\in H,$ $\left\Vert x\right\Vert =1.$

Taking the supremum over $x\in H,$ $\left\Vert x\right\Vert =1,$ we get%
\begin{align*}
\left\Vert T\right\Vert ^{2}& \leq \sup \left\{ \left\vert \left\langle
Tx,x\right\rangle \right\vert ^{2}+2\left\vert \left\langle
Tx,x\right\rangle \right\vert \left( \left\vert \lambda \right\vert -\sqrt{%
\left\vert \lambda \right\vert ^{2}-r^{2}}\right) \right\}  \\
& \leq \sup \left\{ \left\vert \left\langle Tx,x\right\rangle \right\vert
^{2}\right\} +2\left( \left\vert \lambda \right\vert -\sqrt{\left\vert
\lambda \right\vert ^{2}-r^{2}}\right) \sup \left\{ \left\vert \left\langle
Tx,x\right\rangle \right\vert \right\}  \\
& =w^{2}\left( T\right) +2\left( \left\vert \lambda \right\vert -\sqrt{%
\left\vert \lambda \right\vert ^{2}-r^{2}}\right) w\left( T\right) ,
\end{align*}%
which is clearly equivalent to (\ref{4.2}).
\end{proof}

\begin{corollary}
\label{c.5}Let $\varphi ,\phi \in \mathbb{C}$ with $\func{Re}\left( \phi 
\bar{\varphi}\right) >0.$ If $A:H\rightarrow H$ is a bounded linear operator
such that either (\ref{2.6}) or (\ref{2.10}) hold true, then:%
\begin{equation}
\left( 0\leq \right) \left\Vert A\right\Vert ^{2}-w^{2}\left( A\right) \leq 
\left[ \left\vert \phi +\varphi \right\vert -2\sqrt{\func{Re}\left( \phi 
\bar{\varphi}\right) }\right] w\left( A\right) .  \label{4.6}
\end{equation}
\end{corollary}

\begin{remark}
If $M\geq m>0$ are such that either $\left( A^{\ast }-mI\right) \left(
MI-A\right) $ is accreative, or, sufficiently, $\left( A^{\ast }-mI\right)
\left( MI-A\right) $ is self-adjoint and%
\begin{equation}
\left( A^{\ast }-mI\right) \left( MI-A\right) \geq 0\quad \text{in the
operator order,}  \label{4.8}
\end{equation}%
then, by (\ref{3.6}) we have:%
\begin{equation}
\left( 1\leq \right) \frac{\left\Vert A\right\Vert }{w\left( A\right) }\leq 
\frac{M+m}{2\sqrt{mM}},  \label{4.9}
\end{equation}%
which is equivalent to%
\begin{equation}
\left( 0\leq \right) \left\Vert A\right\Vert -w\left( A\right) \leq \frac{%
\left( \sqrt{M}-\sqrt{m}\right) ^{2}}{2\sqrt{mM}}w\left( A\right) ,
\label{4.10}
\end{equation}%
while from (\ref{4.2}) we have%
\begin{equation}
\left( 0\leq \right) \left\Vert A\right\Vert ^{2}-w^{2}\left( A\right) \leq
\left( \sqrt{M}-\sqrt{m}\right) ^{2}w\left( A\right) .  \label{4.11}
\end{equation}%
Also, the inequality (\ref{2.7}) becomes%
\begin{equation}
\left( 0\leq \right) \left\Vert A\right\Vert -w\left( A\right) \leq \frac{1}{%
4}\cdot \frac{\left( M-m\right) ^{2}}{M+m}.  \label{4.12}
\end{equation}
\end{remark}


\begin{thebibliography}{9}
\bibitem{GR} K.E. GUSTAFSON and D.K.M. RAO, \textit{Numerical Range, }%
Springer-Verlag, New York, Inc., 1997.

\bibitem{H} P.R. HALMOS, \textit{A Hilbert Space Problem Book,}
Springer-Verlag, New York, Heidelberg, Berlin, Second edition, 1982.
\end{thebibliography}
\end{document}